\documentclass[12pt,a4paper,english]{article}
\usepackage[T1]{fontenc}
\usepackage[latin1]{inputenc}
\usepackage{babel}
\usepackage{amsthm}
\usepackage{amsmath}
\usepackage{amssymb}
\usepackage{amsfonts}
\usepackage[all]{xy}
\usepackage{hyperref}


\DeclareMathOperator{\Alb}{Alb}

\DeclareMathOperator{\Divf}{\underline{Div}}

\DeclareMathOperator{\HDivf}{\underline{HDiv}}

\DeclareMathOperator{\Expaht}{\mathcal{E} \hspace{-0.5mm} \mathit{xp}}

\DeclareMathOperator{\Frob}{F}
\DeclareMathOperator{\Gal}{Gal}

\DeclareMathOperator{\Hom}{Hom}

\DeclareMathOperator{\Homabk}{Hom_{\Abk}}

\DeclareMathOperator{\Lalb}{Lu}

\DeclareMathOperator{\Nil}{Nil}

\DeclareMathOperator{\Pic}{Pic}

\DeclareMathOperator{\Picf}{\underline{Pic}}

\DeclareMathOperator{\Supp}{Supp}

\DeclareMathOperator{\dv}{div}

\DeclareMathOperator{\id}{id}
\DeclareMathOperator{\im}{im}

\DeclareMathOperator{\modu}{mod}

\DeclareMathOperator{\val}{v}

\newcommand{\Alba}[1]{\Alb_{#1}}

\newcommand{\Albm}[2]{\Alb_{#1,#2}}

\newcommand{\AlbF}[1]{\Alb_{#1,\fmlG}}

\newcommand{\Fm}[2]{\fmlG_{#1,#2}}
\newcommand{\Fmo}[2]{\fmlG_{#1,#2}^{\,0}}

\newcommand{\Hm}[2]{\fmlH_{#1,#2}}

\newcommand{\Lm}[2]{\mathrm{L}_{#1,#2}}
\newcommand{\Lux}[1]{\Lalb_{#1}}

\newcommand{\Luxx}[2]{\Lalb_{#2}^{(#1)}}
\newcommand{\luxx}[2]{\lalb_{#2}^{(#1)}}

\newcommand{\Lalbba}[2]{\Lalb_{#2}^{(#1)}}
\newcommand{\LalbH}[1]{\Lalb_{#1,\fmlH}}
\newcommand{\LalbbH}[2]{\Lalb_{#2,\fmlH}^{(#1)}}

\newcommand{\lalbH}[1]{\lalb_{#1,\fmlH}}
\newcommand{\lalbbH}[2]{\lalb_{#2,\fmlH}^{(#1)}}
\newcommand{\Lalbm}[2]{\Lalb_{#1,#2}}
\newcommand{\lalbm}[2]{\lalb_{#1,#2}}

\newcommand{\Lalbbm}[3]{\Lalb^{(#1)}_{#2,#3}}
\newcommand{\lalbbm}[3]{\lalb^{(#1)}_{#2,#3}}

\newcommand{\Picor}[1]{\Pic^{0,\red}_{#1}}

\newcommand{\HDivfc}[1]{\wh{\HDivf_{#1}}}
\newcommand{\HDivsfc}[2]{\wh{\HDivf_{#1}^{#2}}}

\newcommand{\fundGgeo}[1]{\pi_1^{\mathrm{ab}} \pn{#1}^0}

\newcommand{\fundGgeom}[2]{\pi_1^{\mathrm{ab}} \pn{#1, #2}^0}

\newcommand{\Kaj}[1]{\Sig\pn{#1}}
\newcommand{\NS}[1]{\mathrm{NS}_{#1}}

\newcommand{\Wplus}[2]{\Wpl\pn{#2_1,\ldots,#2_{#1}}}

\newcommand{\Wplusc}[2]{\Wplc\pn{#2_1,\ldots,#2_{#1}}}  
  
\newcommand{\Ab}{\mathsf{Ab}}
\newcommand{\Abk}{\mathcal{A} \mathit{b} / \fld}

\newcommand{\Algk}{\mathsf{Alg} / \fld}

\newcommand{\Artk}{\mathsf{Art} / \fld}

\newcommand{\Mr}{\mathsf{Mr}}

\newcommand{\Mrtaff}[2]{\mathsf{Mr_{\aff}} \pn{#2,#1}}

\newcommand{\Mrmaff}[2]{\Mr_{\aff}\pn{#1,#2}}

\newcommand{\Nat}{\mathbb{N}}
\newcommand{\Zint}{\mathbb{Z}}

\newcommand{\Qrat}{\mathbb{Q}}

\newcommand{\Afn}{\mathbb{A}}
\newcommand{\Prj}{\mathbb{P}}

\newcommand{\Zero}{\mathrm{Z}}

\newcommand{\ab}{\mathrm{ab}}
\newcommand{\aff}{\mathrm{aff}}

\newcommand{\mdl}{D}
\newcommand{\mdll}{E}

\newcommand{\nmr}{j}

\newcommand{\red}{\mathrm{red}}

\newcommand{\lalb}{\mathrm{lu}}

\newcommand{\morphism}{\psi}

\newcommand{\trafo}{\tau}

\newcommand{\clfld}{\overline{\fld}}

\newcommand{\fld}{\mathit{k}}

\newcommand{\fmlG}{\mathcal{F}}

\newcommand{\fmlGr}{\mathcal{G}}
\newcommand{\fmlH}{\mathcal{H}}

\newcommand{\Gm}{\mathbb{G}_{\mathrm{m}}}

\newcommand{\Witt}{\mathrm{W}}
\newcommand{\Wittc}{\wh{\Witt}}
\newcommand{\Wpl}{\Lambda}
\newcommand{\Wplc}{\wh{\Lambda}}

\newcommand{\pntt}{p}

\newcommand{\Ld}{L^{\vee}}

\newcommand{\Xo}{X}

\newcommand{\eps}{\varepsilon}

\newcommand{\lam}{\lambda}

\newcommand{\phe}{\varphi}

\newcommand{\Gam}{\Gamma}

\newcommand{\Sig}{\Sigma}

\newcommand{\fm}{\mathfrak{m}}

\newcommand{\sK}{\mathcal{K}}

\newcommand{\sO}{\mathcal{O}}

\newcommand{\bF}{\mathbb{F}}

\newcommand{\bt}[1]{[#1]}

\newcommand{\pn}[1]{(#1)}
\newcommand{\lrpn}[1]{\left(#1\right)}
\newcommand{\bigpn}[1]{\big(#1\big)}
\newcommand{\Bigpn}[1]{\Big(#1\Big)}
\newcommand{\biggpn}[1]{\bigg(#1\bigg)}

\newcommand{\pair}[1]{\langle #1\rangle}
\newcommand{\lrpair}[1]{\left\langle #1\right\rangle}
\newcommand{\bigpair}[1]{\big\langle #1\big\rangle}

\newcommand{\st}[1]{\{#1\}}
\newcommand{\lrst}[1]{\left\{#1\right\}}

\newcommand{\vrt}[1]{\vert#1\vert}

\newcommand{\ra}{\rightarrow}

\newcommand{\dra}{\dashrightarrow}

\newcommand{\lra}{\longrightarrow}

\newcommand{\Llra}{\Longleftrightarrow}
\newcommand{\lmt}{\longmapsto}

\newcommand{\iso}{\overset{\sim}\longrightarrow}

\newcommand{\lsur}{-\hspace{-.7em} -\hspace{-.7em} \twoheadrightarrow}

\newcommand{\xra}[1]{\xrightarrow{#1}}

\newcommand{\wh}[1]{\widehat{#1}}

\newcommand{\cut}{\cdot}

\newcommand{\tens}{\otimes}

\newcommand{\tms}{\times}

\newcommand{\blank}{\_}
\newcommand{\llul}{\hspace{+.05em} ?}   
\newcommand{\lull}{? \hspace{+.05em}}   

\newcommand{\see}{}
\newcommand{\seecite}{\textrm{see }}
\newcommand{\laurin}{. }                     
\newcommand{\laurink}{, }                   

\newcommand{\vs}{6pt}

\newcommand{\bThm}{\begin{theorem}}
\newcommand{\eThm}{\end{theorem}}
\newcommand{\bAck}{\begin{acknowledgement}}
\newcommand{\eAck}{\end{acknowledgement}}
\newcommand{\bAlg}{\begin{algorithm}}
\newcommand{\eAlg}{\end{algorithm}}
\newcommand{\bAxm}{\begin{axiom}}
\newcommand{\eAxm}{\end{axiom}}
\newcommand{\bCas}{\begin{case}}
\newcommand{\eCas}{\end{case}}
\newcommand{\bClm}{\begin{claim}}
\newcommand{\eClm}{\end{claim}}
\newcommand{\bCcl}{\begin{conclusion}}
\newcommand{\eCcl}{\end{conclusion}}
\newcommand{\bCdn}{\begin{condition}}
\newcommand{\eCdn}{\end{condition}}
\newcommand{\bCjc}{\begin{conjecture}}
\newcommand{\eCjc}{\end{conjecture}}
\newcommand{\bCor}{\begin{corollary}}
\newcommand{\eCor}{\end{corollary}}
\newcommand{\bCrt}{\begin{criterion}}
\newcommand{\eCrt}{\end{criterion}}
\newcommand{\bDef}{\begin{definition}}
\newcommand{\eDef}{\end{definition}}
\newcommand{\bExm}{\begin{example}}
\newcommand{\eExm}{\end{example}}
\newcommand{\bExc}{\begin{exercise}}
\newcommand{\eExc}{\end{exercise}}
\newcommand{\bFct}{\begin{fact}}
\newcommand{\eFct}{\end{fact}}
\newcommand{\bLem}{\begin{lemma}}
\newcommand{\eLem}{\end{lemma}}
\newcommand{\bNot}{\begin{notation}}
\newcommand{\eNot}{\end{notation}}
\newcommand{\bPar}{\begin{para}}
\newcommand{\ePar}{\end{para}}
\newcommand{\bPnt}{\begin{point}}
\newcommand{\ePnt}{\end{point}}
\newcommand{\bPrb}{\begin{problem}}
\newcommand{\ePrb}{\end{problem}}
\newcommand{\bPrp}{\begin{proposition}}
\newcommand{\ePrp}{\end{proposition}}
\newcommand{\bRmk}{\begin{remark}}
\newcommand{\eRmk}{\end{remark}}
\newcommand{\bSol}{\begin{solution}}
\newcommand{\eSol}{\end{solution}}
\newcommand{\bSmr}{\begin{summary}}
\newcommand{\eSmr}{\end{summary}}
\newcommand{\bVar}{\begin{variant}}
\newcommand{\eVar}{\end{variant}}
\newcommand{\bPf }{\begin{prooof}}
\newcommand{\ePf }{\end{prooof}}
\newcommand{\bBf }{\begin{brooof}}
\newcommand{\eBf }{\end{brooof}}

\theoremstyle{plain}
\newtheorem{theorem}{Theorem}[section]
\newtheorem{axiom}[theorem]{Axiom}
\newtheorem{conjecture}[theorem]{Conjecture}
\newtheorem{corollary}[theorem]{Corollary}
\newtheorem{criterion}[theorem]{Criterion}
\newtheorem{lemma}[theorem]{Lemma}
\newtheorem{problem}[theorem]{Problem}
\newtheorem{proposition}[theorem]{Proposition}

\theoremstyle{definition}
\newtheorem{acknowledgement}[theorem]{Acknowledgement}
\newtheorem{algorithm}[theorem]{Algorithm}
\newtheorem{case}[theorem]{Case}
\newtheorem{claim}[theorem]{Claim}
\newtheorem{condition}[theorem]{Condition}
\newtheorem{conclusion}[theorem]{Conclusion}
\newtheorem{definition}[theorem]{Definition}
\newtheorem{example}[theorem]{Example}
\newtheorem{exercise}[theorem]{Exercise}
\newtheorem{fact}[theorem]{Fact}
\newtheorem{notation}[theorem]{Notation}
\newtheorem{point}[theorem]{Point}
\newtheorem{remark}[theorem]{Remark}
\newtheorem{solution}[theorem]{Solution}
\newtheorem{summary}[theorem]{Summary}

\newtheorem{para}[subsubsection]{}

\newenvironment{prooof}[1][Proof]{\textbf{#1.} }
{\ \rule{0.5em}{0.5em}} 
\newenvironment{brooof}[1][Proof]{\textbf{#1.} }
{\ $\square$} 
\newenvironment{variant}[1][Variant]{\textbf{#1.} }
{\ \rule{0.5em}{0.5em}}

\begin{document}


\centerline{ } 
\centerline{\LARGE{Abelianized Fundamental Group of the Affine}} 
\vspace{4pt} 
\centerline{\LARGE{Space over a Finite Field and Big Witt}} 
\vspace{4pt}
\centerline{\LARGE{Vectors in Several Variables}} 
\vspace{25pt} 
\centerline{\large{Henrik Russell}} 
\vspace{10pt} 
\centerline{\large{September 2015}} 
\vspace{25pt} 



\begin{abstract} 
Let $X$ be a normal proper variety over a perfect field $k$. 
We describe abelian coverings of $X$ in terms of 
the functor $\HDivf_X$ of principal relative Cartier divisors on $X$. 
If $k$ is finite, we obtain 
for the maximal abelian extension of the function field of $X$ the relation 
\,$\Gal\lrpn{K^{\ab}_{X} \big| K_{X}\clfld} 
= \Homabk\bigpn{\HDivfc{X},\Gm}$. 
As another application, we present 
the geometric abelianized fundamental group 
of the affine space $\Afn^n$ over a finite field 
by the group of big Witt vectors in $n$ variables, 
a generalization of the (usual) big Witt vectors. 

\end{abstract}


\setcounter{tocdepth}{2}
\tableofcontents{} 

\setcounter{section}{-1}

\newpage 

\section{Introduction} 

Let $X$ be a normal proper variety over a perfect field $\fld$. 
In order to classify abelian coverings of $X$ (in the sense of rational maps), 
Serre observed (basing on Lang's work in the 1950's, see e.g.\ \cite{L_gCFT}) 
that every abelian covering (not arising from an extension of the base field) 
is the pull-back of an isogeny. 
This leads to an explicit description of the class field theory 
of varieties over finite fields, 
using generalized Albanese varieties 
($\seecite$\cite[Thm.~0.4]{Ru2}, \cite[Thm.~3.7]{Ru_CFT}), 
if those are defined (this is the case when $X$ is smooth). 

However, Serre's construction of those isogenies 
($\seecite$\cite[VI, \S~2, No.~8, Prop.~7]{S}) 
involves only \emph{affine} algebraic groups. 
As a consequence, in the class field theory of varieties over finite fields 
one can replace the generalized Albanese varieties 
by their \emph{affinizations}. 
These affinizations are universal objects 
for rational maps from $X$ to commutative affine algebraic groups. 
The structure of these \emph{universal affine groups} 
is much simpler than the one of the generalized Albanese varieties. 
This results in a simplified description 
of class field theory of function fields over finite fields (given in this paper), 
for which we can reduce the assumptions on $X$: 
in particular, smoothness of $X$ can be dropped. 
The price one has to pay for 
getting rid of abelian varieties 
is that we might possibly lose control over the conductor. 

Now suppose the base field $\fld$ is finite, 
let $\clfld$ be an algebraic closure. 
One of the main results of this approach is a 
description of the geometric Galois group 
of the maximal abelian extension of the function field of $X$ 
in terms of the functor $\HDivf_X$ of principal relative Cartier divisors on $X$: 
\[ \Gal\lrpn{K^{\ab}_{X} \big| K_{X}\clfld} 
\;\cong\; \Homabk\bigpn{\HDivfc{X},\Gm} 
\] 
($\see$Theorem \ref{isogeny-of-lux}). 
For any compactification of the affine space $\Afn^n$ 
the generalized Albanese varieties coincide with their affinizations. 
Therefore in this case we still keep control over the conductor 
when applying our method.  
Thus we can use it in order to calculate 
the geometric abelianized fundamental group of $\Afn^n$: 
it is given by the group of $k$-valued points 
of the big Witt vectors in $n$ variables $\Wpl^n$: 
\[ \fundGgeo{\Afn^n} \;\cong\; \Wpl^n\pn{k} 
\] 
(cf.\ Theorem \ref{fundGroup_affSpace}). 
Here $\Wpl^n$ is a generalization of the (usual) big Witt vectors, 
an additive version of which appeared in the context of algebraic K-theory, 
$\seecite$\cite{AGHL}. 
In order to preserve the geometric intuition 
we introduce here an independent multiplicative version: 
the ring of \emph{big Witt vectors in $n$ variables} $\Wpl^n$ 
is the ring-scheme given by the $\Zint$-functor 
\[ R \;\lmt\; 1 + \pn{t_1,\ldots,t_n} \;R[[t_1,\ldots,t_n]] 
\] 
($\see$Definition \ref{biggWitt_def}). 
The affine $\fld$-group $\Wpl^n$ is Cartier dual to its completion at the identity 
$\Wplc^n$ ($\see$Proposition \ref{bigWitt_CartierDual}), 
which is canonically isomorphic to the functor 
\[ R \;\lmt\; \frac{R[t_1,\ldots,t_n]^*}{R^*} 
\] 
($\see$Point \ref{bigWitt-compl}).


\section{Universal Torsors under Affine Groups} 
\label{sec:aff-Univ-Fact-Prbl}

Let $X$ be a proper variety over a perfect field $k$, 
not necessarily smooth or irreducible. 
In \cite[Section 2]{Ru2} we considered 
\emph{categories of rational maps from $X$ 
to torsors  under commutative algebraic groups}. 
In this note we consider 
\emph{categories of rational maps from $X$ 
to torsors  under commutative affine groups}. 
The construction of universal objects 
for the latter categories 
is simpler than for non-affine groups: 
while in the non-affine case we used duality of 1-motives with unipotent part, 
in the case of affine groups it is sufficient to consider 
the Cartier duality between affine groups and formal groups. 
Moreover, in the non-affine case 
the (difficult) group functor $\Divf_X$ of relative Cartier divisors 
was involved, 
while in the case of affine groups we can restrict to 
the (easier) subfunctor $\HDivf_X$ of \emph{principal} relative Cartier divisors.

\subsection{Principal Relative Cartier Divisors} 
\label{relDiv}

The $k$-group functor of principal relative Cartier divisors 
\[ \HDivf_X: \Algk \lra \Ab 
\] 
is the subfunctor of principal elements of 
the $k$-group functor of relative Cartier divisors $\Divf_X$ 
from \cite[No.~2.1]{Ru1}. 
For the purpose of this note it is enough to consider the completion 
\;$\HDivfc{X}: \Artk \lra \Ab$\; of $\HDivf_X$, 
which assigns to a finite $k$-algebra $R$ the group 
\[ \HDivfc{X}\pn{R} 
= \frac{\Gam\bigpn{X \tens R, \pn{\sK_X \tens_k R}^*}}
   {\Gam\bigpn{X \tens R, \pn{\sO_X \tens_k R}^*}} 
= \frac{\pn{\sK_X \tens_k R}^*}{R^*}
\laurink 
\] 
i.e.\ 
\[ \HDivfc{X} 
    = \Gm\pn{\sK_X \tens_k \lull} \big/ \Gm 
    \laurin 
\] 
As this is a quotient of left-exact functors on $\Artk$ 
(= formal $k$-groups), 
and the category of formal $k$-groups is abelian, 
we obtain as in \cite[Prop.~2.1]{Ru2} 

\bPrp 
\label{Divf formal group} 
$\HDivfc{X}$ is a commutative formal $k$-group. 
\ePrp 


\subsection{Rational Maps to Torsors under Affine Groups} 
\label{sub:Categories-of-Rational} 


\bPnt 
\label{induced Trafo}
A rational map $\phe: X \dra L$ to a commutative affine group $L$ 
induces a natural transformation $\trafo_{\phe}: \Ld \lra \HDivfc{X}$: 
let 
\[ \lrpair{\llul,\lull} : \Ld \tms L(\sK_X)/L(k) 
   \lra \Gm(\sK_X\tens\blank)/\Gm 
\] 
be the pairing obtained from Cartier duality, 
then we set 
\[ \trafo_{\phe} = \bigpair{\llul,\bt{\phe}} 
\] 
where $\bt{\phe}$ is the class of $\phe \in L\pn{\sK_X}$ in $L(\sK_X)/L(k)$. 
\ePnt 

\bDef 
\label{Mr_H} 
If $\fmlH$ is a 
formal subgroup of $\HDivf_{X}$, 
denote by $\Mrtaff{\fmlH}{X}$ the category of those rational maps 
for which the image of this induced transformation from Point \ref{induced Trafo} 
lies in $\fmlH$. 
If $\fld$ is an arbitrary perfect base field, 
we define $\Mrtaff{\fmlH}{X}$ via base change to an algebraic closure $\clfld$, 
as in this case we can identify a torsor with the group acting on it. 
\eDef 

A simplification of the proof of \cite[Thm.~0.1]{Ru2} shows 

\bThm 
\label{univ_affObject}
Let $\fmlH$ be a formal (resp.\ dual-algebraic formal) 
$k$-subgroup of $\HDivf_{X}$. 
The category $\Mrtaff{\fmlH}{X}$ admits a universal object 
\[ \lalbbH{1}{X}: X \dra \LalbbH{1}{X} 
\laurin 
\] 
Here $\LalbbH{1}{X}$ is a torsor under an affine (resp.\ algebraic affine) 
commutative $\fld$-group $\LalbbH{0}{X}$, 
which is given by the Cartier dual of \,$\fmlH$. 
\eThm 

\bNot  
In the case $\fmlH = \HDivfc{X}$ 
the category $\Mrtaff{\HDivfc{X}}{X}$ is the category of 
\emph{all} rational maps from $X$ to torsors under affine commutative groups. 
The universal object of this category is denoted by 
\hspace{3mm}$\luxx{1}{X}: X \dra \Luxx{1}{X}$\hspace{3mm} 
(without any specification of $\fmlH$). 
\eNot  

\bRmk 
\label{Alb_constr} 
$\LalbbH{1}{X}$ is generated by $X$ 
and $\LalbbH{0}{X}$ is smooth. 
The rational map 
$\bigpn{\lalbH{X}: X \dra \LalbH{X}} \in \Mrtaff{\fmlH}{X}$ 
is characterized by the fact that the transformation 
\,$\trafo_{\lalbH{X}}: \LalbH{X}^{\vee} \lra \HDivfc{X}$ 
is the identity \,$\fmlH \overset{\id} \lra \fmlH$\, 
(cf.\ \cite[Rmk.~2.18]{Ru2}). 
\eRmk 

\bRmk 
\label{affinization_Alb}
The universal object $\LalbH{X}$ is the affinization 
($\seecite$\cite[III, \S~3, No.~8]{DG}) 
of a generalized Albanese variety $\AlbF{X}$ from \cite{Ru2}, 
if $X$ is smooth, $\fmlH$ is dual-algebraic 
and $\fmlG$ is a dual-algebraic formal subgroup of $\Divf_X$ 
such that $\fmlH = \fmlG \cap \HDivf_X$. 
(This follows directly from the universal property of the affinization.) 
\eRmk 

Via Cartier duality and Galois descent we obtain 
the functoriality of $\Lalbba{i}{X,\fmlH}$ 
(cf.\ \cite[2.3.3]{Ru2}): 

\bPrp 
\label{lalb_H(morphism)}
Let $\fmlH \subset \HDivf_{X}$ be a formal (resp.\ dual-algebraic formal) 
$\fld$-subgroup. 
Let $\morphism: Y \ra X$ be a morphism of proper varieties, 
such that $\morphism\pn{Y}$ meets $\Supp(\fmlH)$ properly. 
Then $\morphism$ induces 
for every formal (resp.\ dual-algebraic formal) $\fld$-subgroup 
$\fmlGr \subset \HDivf_{Y}$ containing $\morphism^*\fmlH$ 
a homomorphism of $\fld$-torsors 
\,$\Lalbba{1}{\morphism,\fmlGr,\fmlH}$ 
and a homomorphism of affine (resp.\ affine algebraic) commutative $\fld$-groups 
\,$\Lalbba{0}{\morphism,\fmlGr,\fmlH}$, 
\[ \Lalbba{i}{\morphism,\fmlGr,\fmlH}: 
   \Lalbba{i}{Y,\fmlGr} \lra \Lalbba{i}{X,\fmlH} 
\hspace{8mm} \textrm{for } i = 1,0 \laurin 
\] 
\ePrp


\subsection{Universal Affine Torsor with Modulus} 
\label{sec:aff-Modulus}

Let $X$ be a proper variety over a perfect field $k$, 
regular in codimension 1. 
Let $\mdl$ 
be an effective Cartier divisor on $X$ 
(possibly non-reduced). 

\bPnt 
\label{filtration}
In \cite{KR2} we assigned to a rational map from $X$ to a torsor 
\,$\phe: X \dra P$\, 
an effective divisor \,$\modu\pn{\phe}$ on $X$, 
the modulus of $\phe$. 
We define a filtration 
\[ \HDivfc{X} = \varinjlim_{\mdl} \Hm{X}{\mdl} 
\] 
(here $\mdl$ ranges over all effective Cartier divisors on $X$) 
by formal subgroups \,$\Hm{X}{\mdl} := \Fm{X}{\mdl} \cap \HDivf_X$\, 
of $\HDivf_{X}$, 
where $\Fm{X}{\mdl}$ are the formal subgroups of $ \Divf_X$ 
from \cite[Def.~3.14]{Ru2}. 
The formal groups $\Hm{X}{\mdl}$ are dual-algebraic 
by \cite[Prop.~3.15 and Lem.~1.17]{Ru2}. 
Moreover they satisfy the following property: 
If $\phe$ maps to a torsor under an affine algebraic commutative group, 
it holds \hspace{4mm}$\modu\pn{\phe} \leq \mdl$\; $\Llra$ 
\;$\im\pn{\trafo_{\phe}} \subset \Hm{X}{\mdl}$, 
see \cite[Lem.~3.16]{Ru2}. 
\ePnt

Using Theorem \ref{univ_affObject} this yields 

\bThm 
\label{AlbMod}
Let $\Mrmaff{X}{\mdl}$ denote the category 
of those rational maps $\phe:X\dra P$ 
to torsors under affine algebraic commutative groups 
s.t.\ \,$\modu\pn{\phe} \leq \mdl$. 
This category admits a universal object 
\[ \lalbbm{1}{X}{\mdl}:X\dra\Lalbbm{1}{X}{\mdl} 
\] 
called the \emph{universal affine torsor of $X$ of modulus $\mdl$}. 
The affine algebraic commutative group $\Lalbbm{0}{X}{\mdl}$ 
acting on $\Lalbbm{1}{X}{\mdl}$ 
is given by the Cartier dual of the formal group $\Hm{X}{\mdl}$. 
\eThm 

Functoriality follows now from Proposition \ref{lalb_H(morphism)}, 
cf.\ \cite[Prop.~3.22 and Cor.~3.23]{Ru2}. 

\bPrp 
\label{lalb_X,D(morphism)} 
Let $\morphism: Y \lra X$ be a morphism of smooth proper varieties. 
Let $\mdl$ be an effective divisor on $X$ 
intersecting $\morphism\pn{Y}$ properly. 
Then $\morphism$ induces 
a homomorphism of torsors $\Lalbba{1}{\morphism,\mdll.\mdl}$ 
and a homomorphism of affine algebraic commutative groups 
$\Lalbba{0}{\morphism,\mdll.\mdl}$, 
\[ \Lalbba{i}{\morphism,\mdll.\mdl}: 
   \Lalbbm{i}{Y}{\mdll} \lra \Lalbbm{i}{X}{\mdl} 
\] 
for each effective divisor $\mdll$ on $Y$ satisfying 
$\mdll \geq \lrpn{\mdl-\mdl_{\red}}\cut Y + \lrpn{\mdl\cut Y}_{\red}$, 
where $B \cut Y$ denotes the pull-back of a Cartier divisor $B$ on $X$ to $Y$. 
\ePrp 

By Cartier duality and Galois descent we have 

\bCor 
\label{Albm(E>D)}
If $\mdl$ and $\mdll$ are effective divisors on $X$ with $\mdll \geq \mdl$, 
then there are canonical surjective homomorphisms 
\;$ \Lalbbm{i}{X}{\mdll} \lsur \Lalbbm{i}{X}{\mdl} 
$\; 
for $i = 1,0$, 
given by $\Lalbba{i}{\id_X,\mdll,\mdl}$. 
\eCor 

\bPrp 
\label{Lux_pro-alg}
The universal object $\Luxx{1}{X}$ 
(for the category of all rational maps to torsors under affine commutative groups) 
is a torsor under a \emph{pro-algebraic} group $\Luxx{0}{X}$. 
\ePrp 

\bPf 
$\Luxx{0}{X}$ is the Cartier dual of 
\,$\HDivfc{X} = \varinjlim \Hm{X}{\mdl}$, 
which is an inductive limit 
of dual-algebraic formal groups $\Hm{X}{\mdl}$ 
($\see$Point \ref{filtration}). 
Thus 
\,$ \Luxx{0}{X} = \varprojlim \Lalbbm{0}{X}{\mdl} 
$\, 
is the projective limit of algebraic affine groups \,$\Lalbbm{0}{X}{\mdl}$. 
\ePf


\section{Geometric Class Field Theory} 
\label{sec:geoCFT}

Let $\fld = \bF_q$ be a finite field, 
$\clfld$ an algebraic closure 
and $X$ an irreducible proper variety over $\fld$. 
We propose to classify abelian coverings of $X$ 
(in the sense of rational maps) 
by principal divisors on $X$. 

\subsection{Geometric Galois Group of a Function Field} 
\label{pull-back_isogeny}

Serre gives in \cite[VI, \S~2, No.~8]{S} 
an explicit construction for any abelian covering 
(arising from a ``geometric situation'') 
as pull-back of an isogeny of \emph{affine groups}: 

\bPnt 
\label{cover-construction}
Let $N$ be a finite group 
and $\morphism: Y \dra X$ a Galois covering (a rational map) 
with Galios group $N$, 
defined and Galois over $\fld$. 
Let $A_N$ be the affine space over $\fld$ defined by the condition that 
$A_N(\clfld) = \clfld[N]$ is the group-algebra of the group $N$. 
Let $G_N$ be the open $\fld$-subvariety of $A_N$, such that 
$G_N(\clfld)$ is the set of invertible elements of $A_N(\clfld)$. 

Then there exists a rational map $\phe: X \dra G_N / N$ over $\fld$ 
such that $Y \dra X$ is isomorphic over $\fld$ 
to the the pull-back of \,$G_N \lra G_N/N$\, via $\phe$. 
(See \cite[VI, \S~2, No.~8, Prop.~7]{S} for a proof.) 
\ePnt 

\bThm 
\label{isogeny-of-lux}
Let $K^{\ab}_{X}$ be the maximal abelian extension 
of the function field $K_{X}$ of $X$. 
The geometric Galois group $\Gal\lrpn{K^{\ab}_{\Xo} \big| K_{X}\clfld}$ 
is isomorphic to the group of $\fld$-valued points of \,$\Lux{X}$, 
which is the Cartier dual of \,$\HDivfc{X}$: 
\[ \Gal\lrpn{K^{\ab}_{X} \big| K_{X}\clfld} 
\;\cong\; \Lux{X}\pn{\fld} 
\;=\; \bigpn{\HDivfc{X}}^{\vee}\pn{\fld} 
\,\laurin 
\] 
\eThm 

\bPf 
Consider a Galois covering $Y \dra X$ as in Point \ref{cover-construction}, 
but suppose that the Galois group $N$ is abelian. 
Then $G_N / N$ is a commutative affine algebraic group, 
hence the rational map $\phe: X \dra G_N / N$ factors through 
$\,\lalbm{X}{\mdl}: X \dra \Lalbm{X}{\mdl}$ 
\footnote{ Here we identify the torsor $\Lalbm{X}{\mdl} := \Lalbbm{1}{X}{\mdl}$ 
with the group $\Lalbbm{0}{X}{\mdl}$ acting on it, 
since $\fld$ is finite and 
hence every $\fld$-torsor admits a $\fld$-rational point.}, 
where $\mdl = \modu\pn{\phe}$. 
Thus $Y \dra X$ is the pull-back of an isogeny onto $\Lalbm{X}{\mdl}$. 
As this isogeny defines an \emph{abelian} extension of function fields, 
it is a quotient of the 
``$q$-power Frobenius minus identity'' $\wp := \Frob_q - \id$ 
(see \cite[VI, No.~6, Prop.~6]{S}). 
If \,$X_{\mdl} \dra X$\, denotes the pull-back of 
\,$\wp: \Lalbm{X}{\mdl} \lra \Lalbm{X}{\mdl}$,  
then $\Gal\pn{Y|X}$ is a quotient of 
\[ \Gal\pn{X_{\mdl}|X} = \ker\pn{\wp} = \Lalbm{X}{\mdl}\pn{\fld} 
\laurin 
\] 
Thus the system $\st{X_{\mdl} \dra X}_{\mdl}$, 
where $\mdl$ ranges over all effective Cartier divisors on $X$, 
is cofinal in the system of geometric abelian Galois coverings 
$\st{Y \dra X}_{Y|X \textrm{ ab.\ geo.}}$. 
As the Galois groups $\Gal\pn{Y|X}$ are finite, 
we have a canonical isomorphism of projective limits 
($\seecite$\cite[Lem.~1.1.9]{RZ}) 
\begin{align*} 
\Gal\lrpn{K^{\ab}_{X} \big| K_{X}\clfld} 
&= \varprojlim_{Y|X \textrm{ ab.\ geo.}} \Gal\pn{Y|X} \\ 
&\cong\, \varprojlim_{\mdl} \Gal\pn{X_{\mdl}|X} 
= \varprojlim_{\mdl} \Lalbm{X}{\mdl}\pn{\fld} 
= \Lux{X}\pn{\fld} 
\laurin 
\end{align*} 
The equality $\Lux{X} = \bigpn{\HDivfc{X}}^{\vee}$ 
is due to the construction of $\Lux{X}$ ($\see$Thm.~\ref{univ_affObject}). 

\hfill 
\ePf

\subsection{Fundamental Group of the Affine Space} 
\label{Exm-affSp}

Finally we want to apply the tools developed so far 
to the class field theory of some given concrete variety. 
The affine space -- 
a priori the easiest example one can imagine -- 
is a welcome example 
in order to test our theory. 

Let $\fld = \bF_q$ be a finite field. 
The abelianized fundamental group of the affine line $\Afn^1$ over $\fld$ 
is supposed to be well-known 
since the generalized Jacobian of Rosenlicht is established. 
(Strangely I did not find this example in any text-book.) 
In this subsection we want to show that 
the abelianized geometric fundamental group 
of the affine space $\Afn^n$ over $\fld$ 
is just as easy as the case $n = 1$. 

\bThm 
\label{fundGroup_affSpace}
Let $\Afn^n$ be the $n$-dimensional affine space over $\fld = \bF_q$. 
Let $t_1,\ldots, t_n$ be affine coordinates of $\Afn^n$. 
Then 
\[ \fundGgeo{\Afn^n} 
= \Bigpn{1 + \pn{t_1^{-1},\ldots, t_n^{-1}} \,k[[t_1^{-1},\ldots, t_n^{-1}]]}^{\tms} 
\,\laurin 
\] 
In other words: the abelianized geometric fundamental group 
of $\Afn^n$ over $\bF_q$ 
is isomorphic to the group of $\bF_q$-valued points 
of the big Witt vectors in $n$ variables 
(see Section \ref{sec:bigWitt-sevVar}). 
\eThm 

\bPf 
Consider $\Afn^n$ as an open subvariety of $X := \pn{\Prj^1}^n$, 
the product of $n$ copies of the projective line $\Prj^1$. 
Let $S \subset X$ be the divisor at infinity 
such that $\Afn^n = X \setminus S$. 
Since the N\'eron-Severi group $\NS{\Prj^1}$ of $\Prj^1$ is torsion-free, 
the same is true for $\NS{X}$ and it holds 
\[ \Kaj{X} 
:= \Hom\Bigpn{\varinjlim_{n} \Homabk\pn{\mu_n, \NS{X}}, \Qrat/\Zint} 
= 0 
\,\laurink 
\] 
hence for every effective divisor $\mdl$ on $X$ we have 
\[ \fundGgeom{X}{\mdl} = \Albm{X}{\mdl}\pn{\fld} 
\] 
by \cite[Cor.~3.14]{Ru_CFT}.\footnote{ 
For technical reasons 
the statement is formulated there only for $\dim X \leq 2$ 
but holds for arbitrary dimension. 
In the case $X = \pn{\Prj^1}^n$ and 
$\vrt{\mdl} = \bigcup_{i = 1}^{n} \Prj^1 \tms \ldots \tms \st{\infty} \tms \ldots \tms \Prj^1$ 
the necessary induction argument is particularly easy.} 
Since \;$\Alba{X} = \pn{\Picor{X}}^{\vee} = 0$\; 
we have for every $\mdl$ 
\[ \Albm{X}{\mdl} = \Lalbm{X}{\mdl} 
\laurin 
\] 
Then 
\begin{align*} 
\fundGgeo{X \setminus S} 
&= \varprojlim_{\substack{\mdl \\ \vrt{\mdl} \subset S}} \fundGgeom{X}{\mdl} 
= \varprojlim_{\substack{\mdl \\ \vrt{\mdl} \subset S}} \Albm{X}{\mdl}\pn{\fld} \\ 
&= \varprojlim_{\substack{\mdl \\ \vrt{\mdl} \subset S}} \Lalbm{X}{\mdl}\pn{\fld} 
= \Bigpn{\varinjlim_{\substack{\mdl \\ \vrt{\mdl} \subset S}} 
      \Hm{X}{\mdl}}^{\vee} \pn{\fld} \\ 
&= \Bigpn{\HDivsfc{X}{S}}^{\vee} \pn{\fld} 
\end{align*} 
where $\HDivf_{X}^S$ denotes the functor of principal relative Cartier divisors 
with support in $S$. 
Then $\HDivsfc{X}{S}$ is the functor 
that assigns to a finite $\fld$-ring $R$ the group 
\begin{align*} 
\HDivsfc{X}{S}\pn{R} 
&= \frac{\Gam\pn{\sO_{\Afn^n} \tens_k R}^*}{R^*} 
= \frac{R[t_1,\ldots,t_n]^*}{R^*} 
\,\laurin 
\end{align*} 
Thus, as we will see in 
Point \ref{bigWitt-compl}, 
\[ \HDivsfc{X}{S} = \Wplusc{n}{t} 
\] 
is the completion at $1$ of the $\fld$-group $\Wplus{n}{t}$ 
of big Witt vectors in $n$ variables, 
and this completion is the Cartier dual of $\Wplus{n}{t^{-1}}$ 
by Prop.~\ref{bigWittn_geoDual}. 
Then 
\begin{align*} 
\Bigpn{\HDivsfc{X}{S}}^{\vee} \pn{\fld} 
&= \Wplus{n}{t^{-1}} \pn{\fld} \\ 
&= \Bigpn{1 + \pn{t_1^{-1},\ldots, t_n^{-1}} \,k[[t_1^{-1},\ldots, t_n^{-1}]]}^{\tms} 
\laurin 
\end{align*} 

\hfill 
\ePf

\section{Big Witt Vectors in Several Variables} 
\label{sec:bigWitt-sevVar}

This section is devoted to a generalization of the ring-scheme of 
big Witt vectors (in one variable). 
Big Witt vectors in several variables 
can be expressed as an infinite product of 
big Witt vectors in one variable. 
We compute the Cartier dual of this group-scheme 
via reduction to the one-variable case. 

Let $\fld$ be a ring. 

\bDef 
\label{biggWitt_def}
The group of \emph{big Witt vectors in $n$ variables} 
is the $\fld$-group defined by 
\[ \Wpl^n := \Wplus{n}{t} 
:= \ker\biggpn{\Gm\Bigpn{\llul \tens \fld[[t_1,\ldots,t_n]]} \lra \Gm} 
\] 
w.r.t.\ the augmentation map 
\[ \fld[[t_1,\ldots,t_n]] \lra \fld \,\laurink 
\hspace{12mm} t_i \lmt 0 \,\laurink 
\] 
i.e.\ $\Wpl^n$ is the functor that assigns to a $\fld$-algebra $R$ the group 
\[ \Wpl^n\pn{R} = \Bigpn{1 + \pn{t_1,\ldots,t_n} \,R[[t_1,\ldots,t_n]] }^{\tms} 
\,\laurin 
\] 

Every element \,$\lam \in \Wpl^n\pn{R}$\, has a unique decomposition 
\[ \lam = \prod_{\nu_1,\ldots,\nu_n} 
\bigpn{1 - r_{\nu_1,\ldots,\nu_n}\, t_1^{\nu_1} \cdots t_n^{\nu_n}} 
\,\laurin 
\] 
\eDef 

\bRmk 
The case $n = 1$ is the $\fld$-ring of (usual) big Witt vectors $\Wpl := \Wpl^1$. 
The ring structure is given as follows: 
addition on $\Wpl$ is the multiplication of formal power series, 
while multiplication on $\Wpl$ is given by 
\begin{eqnarray*}
\Wpl \tms \Wpl & \overset{*} \lra & \Wpl \\ 
\prod_{i \geq 1} \pn{1 - a_i t^i} \,,\, \prod_{j \geq 1} \pn{1 - b_j t^j} & \lmt & 
\prod_{i,j \geq 1} \Bigpn{1 - a_i^{\frac{j}{(i,j)}} b_j^{\frac{i}{(i,j)}} t^{\frac{ij}{(i,j)}} }^{(i,j)} 
\end{eqnarray*} 
where $(i,j) := \gcd(i,j)$, 
cf.\ \cite[I, \S~1, Prop.~(1.1)]{Bl1}. 
\eRmk 

\bPrp 
\label{bigWitt_decomposition}
The group of big Witt vectors $\Wpl^n$ is canonically isomorphic to
 an (infinite) product of copies of $\Wpl$: 
\[ \Wpl^n \;\iso 
\prod_{\substack{\nu_1,\ldots,\nu_n \\ \gcd(\nu_1,\ldots,\nu_n) = 1}} \Wpl 
\] 
\ePrp 

\bPf 
Let $\nu := \pn{\nu_1,\ldots,\nu_n}$ be a multi-index, 
write $t^{\nu} := t_1^{\nu_1} \cdots t_n^{\nu_n}$. 
Then we obtain a canonical isomorphism 
\begin{eqnarray*} 
\Wplus{n}{t} & \cong & \prod_{\substack{\nu \\ \gcd(\nu)=1}} \Wpl\pn{t^{\nu}} \\ 
\prod_{\nu} \bigpn{1 - r_{\nu} \,t^{\nu}} 
& = & \prod_{\substack{\nu \\ \gcd(\nu)=1}} \prod_{i} \bigpn{1 - r_{i \nu} \,t^{i \nu}} 
\end{eqnarray*} 

\hfill 
\ePf 

\bPrp 
\label{bigWitt_properties}
The functor $\Wpl^n$ is a smooth connected unipotent commutative $\fld$-group 
and a $\fld$-ring. 
\ePrp 

\bPf 
Due to the decomposition of $\Wpl^n$ from Prop.~\ref{bigWitt_decomposition}, 
this follows from the properties of the (usual) big Witt vectors $\Wpl$. 
\ePf 

\vspace{\vs} 
Now let $\fld$ be a field. 

\bPnt 
\label{bigWitt-compl}
The \emph{completion at 1} of the big Witt vectors in $n$ variables 
is the formal $\fld$-subgroup of $\Wpl^n$ 
\[ \Wplc^n \,:=\, \Wplusc{n}{t} \;\subset\; \Wplus{n}{t} 
\] 
that assigns to a $\fld$-algebra $R$ the group 
\[ \Wplc^n\pn{R} = 
\lrst{ 1 + \hspace{-2mm} \sum_{\nu_1,\ldots,\nu_n} 
    r_{\nu_1,\ldots,\nu_n}\, t_1^{\nu_1} \cdots t_n^{\nu_n} 
\left| 
\begin{array}{l} 
r_{\nu_1,\ldots,\nu_n} \in \Nil\pn{R} \;\textrm{ for all } \;{\nu_1,\ldots,\nu_n} \\ 
r_{\nu_1,\ldots,\nu_n} = 0 \;\textrm{ for almost all } \;{\nu_1,\ldots,\nu_n} 
\end{array} 
\right. 
} 
\laurin 
\] 
Then \,$\Wplc^n$\, is canonically isomorphic to the $\fld$-group functor 
\[ R \;\lmt\; \frac{R[t_1,\ldots,t_n]^*}{R^*} 
\] 
i.e.\ 
\[ \Wplusc{n}{t} \;\cong\; \frac{\Gm\bigpn{\llul \tens \fld[t_1,\ldots,t_n]}}{\Gm} 
\;\laurin 
\] 
\ePnt 

\bPrp 
\label{bigWitt_CartierDual}
The big Witt vectors $\Wpl^n$ are Cartier dual to $\Wplc^n$, 
its completion at 1: 
\[ \bigpn{\Wpl^n}^{\vee} \;\cong\; \Wplc^n 
\laurin 
\]  
The pairing is given by the multiplication $*$ on $\Wpl^n$ 
composed with an evaluation at \,$t = 1$: 
\[ 
\begin{array}{ccccl} 
\Wplc^n \tms \Wpl^n & \overset{*} \lra & \Wplc^n & \xra{\;t = 1\;} & 
\hspace{3.5mm} \Gm \\ 
f \;,\; g & \lmt & f * g & \lmt & \pn{f * g}(1) \laurin 
\end{array} 
\] 
\ePrp 

\bPf 
According to the decomposition of $\Wpl^n$ 
from Prop.~\ref{bigWitt_decomposition}, 
we can reduce to the case of (usual) big Witt vectors $\Wpl$. 
We have an isomorphism 
\[ 
   \Pi\eps: 
   \prod_{\substack{\nmr \geq 1 \\ \pn{\nmr,p}=1}} \Witt \iso \Wpl \;, 
   \hspace{12mm} 
   \pn{v_{\nmr}}_{\nmr} \lmt \prod_{\nmr} \Expaht\bigpn{v_{\nmr}, t^{\nmr}} 
\] 
where $\Expaht$ denotes the Artin-Hasse exponential, 
from an infinite product of copies of the ring of (small) Witt vectors $\Witt$ 
to the ring of big Witt vectors $\Wpl$ 
($\seecite$\cite[III, No.~1, Prop.\ on p.~53]{D}). 
The Cartier dual of $\Witt$ is given by $\Wittc$, its completion at 0, 
and the pairing is 
\[ 
\begin{array}{ccccccl} 
\Wittc \tms \Witt & \overset{\cdot} \lra & \Wittc & \xra{\;\Expaht\;} & \Wplc & 
\xra{\;t = 1\;} & \hspace{3.5mm} \Gm \\ 
v \;,\; w & \lmt & v \cdot w & \lmt & \Expaht\pn{v \cdot w, t} & \lmt & 
\Expaht\pn{v \cdot w, 1} 
\end{array} 
\] 
($\seecite$\cite[V, \S~4, Cor.~4.6]{DG}). 
As the map \,$\Pi\eps$\, is a homomorphism of $\fld$-rings 
($\seecite$\cite[V, \S~5, Thm.~5.5]{DG}), 
the result follows. 
\ePf 

\bPrp[geometric description of Cartier duality] 
\label{bigWitt-duality_geom-Interpretation}
Let $\fld$ be a perfect field, $\clfld$ an algebraic closure. 
The Cartier duality between $\Wplc$ and $\Wpl$ 
is expressed by the pairing 
\;\;$\Wplc(t) \tms \Wpl(t^{-1}) \lra \Gm$\;\; 
determined by 
\begin{eqnarray*} 
\Wplc(t)(R) \tms \Wpl(t^{-1})(\clfld) & \lra & \Gm(R \tens \clfld) \\ 
f \;,\; g \hspace{15mm} & \lmt & \prod_{\pntt \in \vrt{\Zero(g')}} f(\pntt)^{\val_{\pntt}(f)} 
\end{eqnarray*} 
where \,$g' \in 1 + t^{-1} \clfld[t^{-1}]$\, is a truncation of \;$g \mod t^{-m} \clfld[[t^{-1}]]$ 
for $m$ sufficiently large, 
and \,$\Zero(g')$\, is the divisor of zeroes of \,$g'$. 
\ePrp 

\bPf 
Consider the Jacobian of $C = \Prj^1$ over $\clfld$ 
with modulus $\mdl = n \,\bt{\infty}$.
Its affine part is 
\[ \Lm{\Prj^1}{n \bt{\infty}} \pn{\clfld} 
 = \frac{\sO_{\Prj^1,\infty}^*}{\fld^* \tms \bigpn{1+\fm_{\Prj^1,\infty}^{n}}} 
 = \frac{1+\fm_{\Prj^1,\infty}}{1+\fm_{\Prj^1,\infty}^{n}} 
 = \frac{1 + t^{-1} \,\clfld[[t^{-1}]]}{1 + t^{-n} \,\clfld[[t^{-1}]]} 
\,\laurin 
\]             
Moreover, as $\Picf_{\Prj^1} = 0$, 
we have $\Divf^0_{\Prj^1} = \HDivf_{\Prj^1}$ 
and $\Fmo{\Prj^1}{n \bt{\infty}} = \Hm{\Prj^1}{n \bt{\infty}}$. 
Then 
\begin{align*} 
& \Wplc(t)\pn{R} \cong \frac{R[t]^*}{R^*} = \HDivsfc{\Prj^1}{\infty} \pn{R} 
= \varinjlim_n \Fmo{\Prj^1}{n \bt{\infty}}\pn{R} 
\,\laurink \\  
& \Wpl(t^{-1})\pn{\clfld} = 1 + t^{-1} \,\clfld[[t^{-1}]] 
= \varprojlim_n \Lm{\Prj^1}{n \bt{\infty}} \pn{\clfld} 
\,\laurin 
\end{align*} 
According to \cite[Prop.~2.5, Point 2.7 and Point 2.4]{Ru_CFT} 
the pairing \\ 
\,$\Wplc(t) \tms \Wpl(t^{-1}) \lra \Gm$\hspace{3mm} 
is thus given by the local symbol at infinity 
\[ \pair{f,g}_{\Prj^1,\infty} = \pn{f,g}_{\infty}^{-1} = f \bigpn{\dv(g')} 
\,\laurin 
\footnote{ 
The pairing $\pair{\llul,\lull}_{\Prj^1,\infty}$ here 
differs from the pairing considered in \cite{Ru_CFT} 
by an inversion in $\Gm$.} 
\] 
Here we can truncate $g$ as stated, 
since $f$ annihilates $m^{\rm th}$ higher units 
if \,$f \in \Fmo{\Prj^1}{m \bt{\infty}}$ for some $m \in \Nat$ 
($\seecite$\cite[\S~6, Prop.~6.4 (3)]{KR2}). 
Now the poles of \,$g \in 1 + t^{-1} \,\clfld[[t^{-1}]]$\, lie all at $0$, 
and for \,$f \in 1 + t \,\clfld[t]$\, it holds \,$f(0) = 1$, 
hence 
\[ f \bigpn{\dv(g')} = f \bigpn{\Zero(g')} 
= \prod_{\pntt \in \vrt{\Zero(g')}} f(\pntt)^{\val_{\pntt}(f)} 
\;\laurin 
\] 

\hfill 
\ePf 

\vspace{\vs} 

Combining the geometric description 
from Prop.~\ref{bigWitt-duality_geom-Interpretation} 
and the decomposition from Prop.~\ref{bigWitt_decomposition} 
we obtain 

\bPrp 
\label{bigWittn_geoDual}
Over a perfect field $\fld$ we have the following duality: 
\[ {\Wplus{n}{t^{-1}}}^{\vee} \;=\; \Wplusc{n}{t} 
\] 
and the pairing between \,$\Wplus{n}{t^{-1}}$\, and \,$\Wplusc{n}{t}$\, 
is induced by the local symbol map. 
\ePrp

\textbf{Acknowledgement.} 
The the roots of this paper originate from my fellowship time in Kyoto. 
I thank Kazuya Kato for his hospitality and helpful discussions. 
It was Kato's idea to erase abelian varieties from the theory 
of generalized Albanese varieties.

\vspace{1mm} 

\begin{flushright} 
e-mail: \texttt{henrik.russell@gmail.com} 
\end{flushright}

\end{document}